\documentclass[letterpaper,11pt]{article}

\usepackage[caption=false]{subfig}
\usepackage[dvips]{graphicx}
\graphicspath{{figures/}}
\usepackage{amsmath, amssymb}
\usepackage{tikz} \usetikzlibrary{calc}
\usepackage{psfrag}
\usepackage{pgfplots}
\usepackage{algorithm}     
\usepackage{algpseudocode} 
\usepackage{arydshln}
\usepackage{listings}

\DeclareMathOperator*{\diag}{diag}

\DeclareMathOperator*{\blkdiag}{blkdiag}

\title{Automated Model Generation for Analysis of Large-scale Interconnected Uncertain Systems}
\author{F. Karami, S. Khoshfetrat Pakazad, A. Hansson and A. Afshar}

\begin{document}
\maketitle
\begin{abstract}
The first challenge in robustness analysis of large-scale interconnected uncertain systems is to provide a model of such systems in a standard-form that is required within different analysis frameworks. This becomes particularly important for large-scale systems, as analysis tools that can handle such systems heavily rely on the special structure within such model descriptions. We here propose an automated framework for providing such models of large-scale interconnected uncertain systems that are used in Integral Quadratic Constraint (IQC) analysis. Specifically, in this paper we put forth a methodological way to provide such models from a block-diagram and nested description of interconnected uncertain systems. We describe the details of this automated framework using an example.
\end{abstract}

\section{Introduction}
\label{sec:int}
Over the past few years, there has been remarkable research interest in analysis of interconnected uncertain systems, see for example \cite{IEEEhowto:Langbort}, \cite{IEEEhowto:Andrea}, \cite{IEEEhowto:Dullerud}, \cite{IEEEhowto:Langbort1}, \cite{IEEEhowto:Langbort2}, \cite{IEEEhowto:Ugrinovskii}, \cite{IEEEhowto:Fang}, \cite{IEEEhowto:Motee}, \cite{IEEEhowto:Chen} and \cite{IEEEhowto:Scorletti}. Robustness analysis of large-scale interconnected uncertain systems poses different  computational challenges mainly due to their large state, uncertainty and input-output dimensions. There are different approaches for analyzing robustness of interconnected uncertain systems, among which, IQC analysis is one of most general ones which provides a unified framework for analyzing uncertain systems with different types of uncertainties,\cite{IEEEhowto:Megretski1}, \cite{IEEEhowto:jonsson} and \cite{IEEEhowto:Zhou}. As an example of large-scale interconnected systems, one can refer to power grids. These networks can be viewed as an interconnection of of huge number of subsystems or components. Due to the continuing deregulation in generation of power and distribution sectors in the U.S. and Europe, the number of players involved in the generation and distribution of power has increased significantly. This can in turn the introduce uncertainties within the system, e.g., in generation and consumption. Due to this intrinsic uncertainty their robustness analysis is put in new demand.

In\cite{IEEEhowto:Sina1}, \cite{IEEEhowto:Sina2} the authors proposed algorithms to alleviate computational issues with conventional robustness analysis techniques when applied to large-scale interconnected uncertain systems. Due to a common structure in such systems, that is sparsity of interconnections among subsystems, it is possible to devise efficient centralized or distributed algorithms for solving the corresponding analysis problems, see e.g., \cite{IEEEhowto:Sina2},\cite{IEEEhowto:Demourant}. This is done by exploiting the sparsity in the interconnections which then makes it possible to either use efficient sparse solvers or decompose the problem and use distributed computational methods \cite{IEEEhowto:Sina2,IEEEhowto:Sina3}. Such efficient analysis algorithms rely on a particular model description. Consequently, for analyzing large-scale interconnected uncertain systems with sparse interconnections, it is necessary to provide the algorithms with such model descriptions. However, generally the provided models are not in the suitable format and are expressed using nested block-diagram-based representations within commonly used simulation softwares, such as Simulink and Modelica. This description allows for natural modelling of such systems. Considering the sheer size of such models then, it is essential to devise automated approaches for extracting the model required for the analysis algorithms from such descriptions.

In this paper, we first put forth a general description of interconnected uncertain systems. However initially this description differs from the standard representation of uncertain systems, i.e., LFT representation, and hence we are prohibited from using efficient analysis tools discussed in \cite{IEEEhowto:Sina3} for analyzing such systems. In order to solve this issue, we show how to convert the aforementioned mathematical description to that of used in \cite{IEEEhowto:Sina2}. To this end, we have taken extra care to preserve the sparsity in the interconnections in the modified description. In order to make the conversion to the standard form, we present a two stage procedure. This procedure is also illustrated using an example.

\section{Mathematical Description of Interconnected Uncertain Systems}
\label{sec:2}
\begin{figure}
\centering
 \includegraphics[width=12cm]{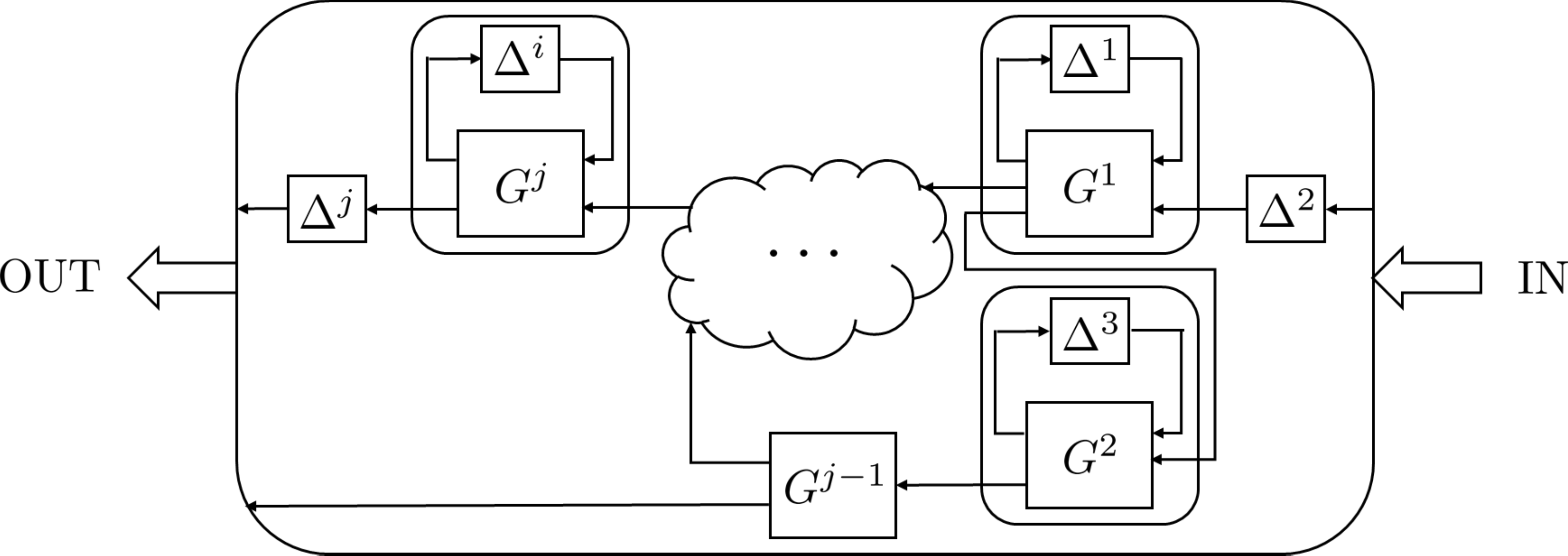}
 \caption{ The structure of uncertain interconnected system}\label{fig:2}
\end{figure}

In this section, we put forth a mathematical description for networks of dynamical systems in the presence of uncertainty. Figure  \ref{fig:2} illustrates one such system. As shown in Figure \ref{fig:2}, main building blocks for describing such systems are dynamic interconnected systems $G^i$, and uncertainty blocks $\Delta^j$. Generally, the block-diagram-based description of interconnected uncertain systems can also include nested structure. That is certain blocks are themselves indirectly described using these fundamental blocks. Due to this, the description of the total system in its most basic form also relies on these two types of blocks. First we provide a formal description of these fundamental blocks.

Let us assume that the interconnected system includes $N$ dynamic system blocks. A dynamic system block is defined using a system transfer function matrix $G^i \in \mathcal{RH}_{\infty}^{p_i \times m_i}$ with $u^i$ and $y^i$, as its $m_i$-dimensional and $p_i$- dimensional input and output vectors, respectively, and is given as
\begin{equation} \label{EQ__1_}
y^i=G^i u^i,
\end{equation}
for $i = 1,\dots, N$. The uncertainty blocks represent causal, bounded operators, $\Delta^i : R^{d_i}\to R^{d_i}$, with $d_i$-dimensional input and output vectors $p^i$ and $q^i$, respectively, and is given as
\begin{equation}\label{EQ__2_}
q^i = \Delta^i( p^i),
\end{equation}
for $i = 1,\dots, d$. Here $d$ is the number of such blocks in the interconnected uncertain system representation. The uncertainty blocks are commonly described using IQCs, \cite{IEEEhowto:Megretski1,IEEEhowto:jonsson}. 
Let us now summarize all the input-output relationship of the dynamic system and uncertainty blocks of interconnected system as
\begin{subequations} \label{EQ__5ab_}
\begin{align}
y&=Gu,\\
q&=\Delta(p),
\end{align}
\end{subequations}
where $u = (u^1, \dots, u^N)$, $y = (y^1, \dots, y^N)$, $q = (q^1, \dots, q^d)$ and $p = (p^1, \dots, p^d)$. Moreover, $G$ and $\Delta$ are block-diagonal matrices, with $G^i$ and $\Delta^j$ as their diagonal elements, respectively, for $i=1, 2,\dots, N$ and $j=1, 2, \dots, d$. This means that $G\in RH_{\infty} ^ {\bar{p}\times\bar{m}}$, with $\bar{m}=\sum _{i=1}^{N} m_i , \bar p=\sum _{j=1}^{d} p_i $, and $\Delta$ :$ R^{\bar d} \to R^{\bar d} $ with $\bar d=\sum _{i=1}^{d} d_i $. What remains now is to describe the interconnection matrix among the aforementioned fundamental blocks.

Having described the building blocks of the system, we now discuss the interconnection among these blocks. This can be done using a so-called interconnection constraint. To this end, we use a 0-1 matrix $\Gamma$ as below
\begin{equation} \label{EQ_6_}
\tilde u=\tilde \Gamma \tilde y
\end{equation}
where $ \tilde u = (u, p)$ and $ \tilde y = (y, q)$ and $\tilde \Gamma$ is a $0-1$ matrix that. The structure of $\tilde \Gamma$ can be described in more detail as
\begin{align}\label{EQ_7_}
\begin{bmatrix} u \\ p \end{bmatrix} = \begin{bmatrix}  \Gamma_{gg} & \Gamma_{gd} \\  \Gamma_{dg} & \Gamma_{dd} \end{bmatrix} \begin{bmatrix} y \\ q \end{bmatrix}
\end{align}
and it basically describes how outputs of each of the blocks are connected to inputs of other blocks. Given $\tilde \Gamma$ or the description of interconnections together with \eqref{EQ__5ab_} completes our description of interconnected uncertain systems. It is also possible to accommodate external inputs and outputs to and from an interconnected uncertain system using a similar approach. Let us define $\bar u=(u, p, u_{\textrm{ext}})$ and $\bar y=(y, q, y_{\textrm{ext}})$, where $u_{\textrm{ext}}$ and $y_{\textrm{ext}}$ are vectors of the external inputs and outputs, respectively. Then we can incorporate the external inputs and outputs by modifying the interconnection description as
\begin{equation} \label{EQ_8_}
\bar u= \bar \Gamma \bar y
\end{equation}
where
\begin{align}\label{EQ_9_}
\bar \Gamma = \begin{bmatrix} \Gamma_{gg}& \Gamma_{gd}& \Gamma_{gi}\\ \Gamma_{dg} & \Gamma_{dd} & \Gamma_{di}\\  \Gamma_{og} & \Gamma_{od} &  \Gamma_{oi} \end{bmatrix}.
\end{align}
With the description of the interconnections in place we can now summarize the interconnected uncertain systems description as
\begin{subequations} \label{EQ_10abc_}
\begin{align}
y&=Gu,\\
q&=\Delta p,\\
\bar u &=\bar \Gamma \bar y.
\end{align}
\end{subequations}
Next we discuss how this description of interconnected uncertain systems in Figure \ref{fig:5} can be transformed to a linear fractional representation used within scalable analysis algorithm.
\section{System Description Conversion to LFT Rrepresentation}
\label{sec:sec-2-3}
\begin{figure}
        \centering
        \includegraphics[width=2.5cm]{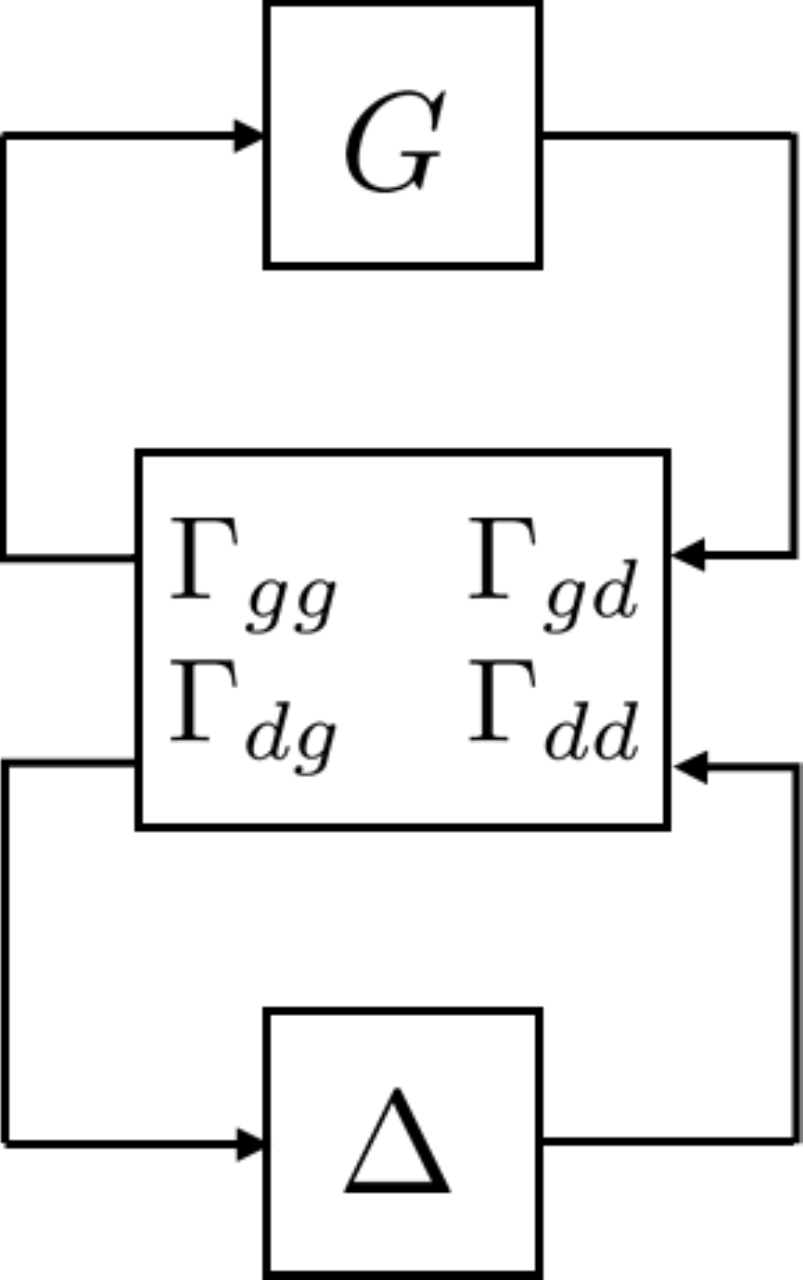}
        \caption{The structure of the system}
        \label{fig:5}
 \end{figure}

\begin{figure}
        \centering
        \includegraphics[width=2.5cm]{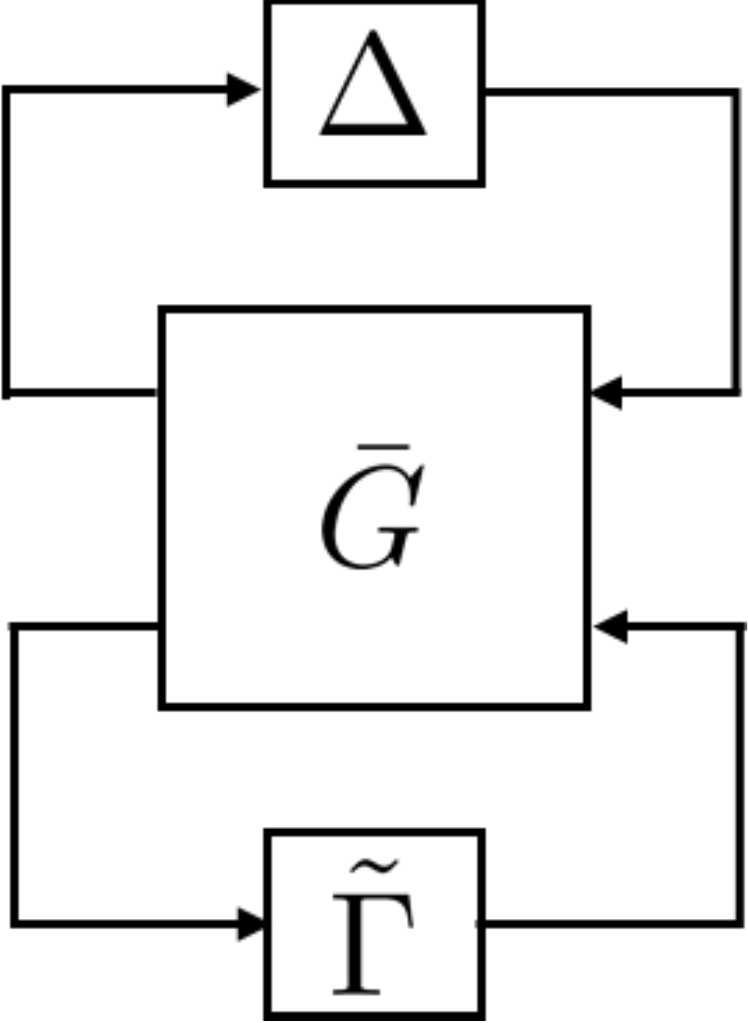}
        \caption{Final structure of the system}
        \label{fig:6}
\end{figure}

At this point we have defined a general mathematical description for interconnected uncertain systems. For the sake of simplicity and brevity we assume that we do not have any external inputs and outputs. This is particularly the case when we are concerned with analysis only. In this section we describe how to convert the description presented above, also illustrated as in Figure \ref{fig:5}, to a standard LFT representation. To this end, let us first review this standard model description for interconnected uncertain systems. The system description considered in scalable analysis algorithms are given in the form
\begin{equation}\label{EQ__4abc_}
\begin{split}
p &= G_{pq} q + G_{pw} w \\
z &= G_{zq} q + G_{zw} w \\
q &= \Delta( p)\\
w &= \tilde \Gamma z,
\end{split}
\end{equation}
e.g., see \cite{IEEEhowto:Sina1}. Now let us discuss how the representation discussed in the previous section can be written in the same format as above. Let us first define
\begin{equation}\label{EQ_11_}
\bar G = \begin{bmatrix} 0 & 0 & I \\ 0 & G & 0 \\ I & 0 & 0 \end{bmatrix}\\
\end{equation}
Then it is possible to rewrite the interconnected system description in Section \ref{sec:2} as in \eqref{EQ__4abc_}, also shown in Figure \ref{fig:6}. Moreover, notice that $\bar G \in \mathcal{RH}_\infty$.
It is then possible to identify different elements in \eqref{EQ__4abc_} as
\begin{equation}
\begin{split}
G_{pq} &= 0, \quad \quad \  G_{pw} = \begin{bmatrix} 0 & I \end{bmatrix}\\
G_{zq} &= \begin{bmatrix} 0 \\ I \end{bmatrix}, \quad G_{zw} = \begin{bmatrix}G & 0 \\ 0 & 0 \end{bmatrix}
\end{split}
\end{equation}
Consequently this description is now in LFT format. So far, we have described how given a general mathematical description of an interconnected uncertain system, it is possible to convert this description to a standard LFT representation. However, this mathematical description, particularly the interconnection description, is not generally given and it needs to be extracted from a nested block-diagram schematic of such systems. For large-scale interconnected uncertain systems, this extraction cannot be done manually and needs to be automated. Next we describe our proposed automated approach to extract this mathematical description.

\section{Nested Description of Interconnected Uncertain Systems and Interconnection Matrix Extraction}\label{sec:3}

\begin{figure}
        \centering
        \includegraphics[width=9cm]{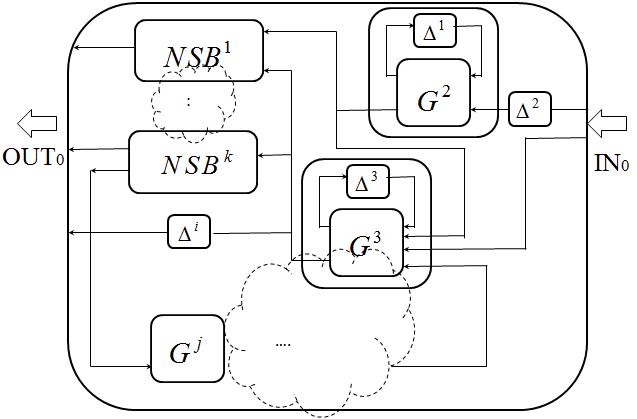}
        \caption{General structure of an interconnected large-scale system}
        \label{fig:7}
\end{figure}

In order to generate a model for a large-scale interconnected uncertain system in the from of \eqref{EQ_10abc_}, we need to extract the interconnection matrix $\bar \Gamma$ from a nested description of the system, as shown in Figure \ref{fig:7}. In order to be able to collapse the nested structure and recover the description in (\ref{EQ_10abc_}), we first need to describe this structure in the system. This can commonly be done using a graphical representation, particularly a tree. Here we assume that this representation is given. This can be provided using different approaches such as regional or adaptive clustering \cite{IEEEhowto:Ragusa}. Next, we describe this graphical description, and then we discuss how it can be used to recover the model as in \eqref{EQ_10abc_}.
\subsection{A Graphical Presentation of the Nested Structure}
\label{sec:sub-3-1}
As was mentioned before the provided models for interconnected uncertain systems have a nested structure. This is mainly to make descriptions of interconnected uncertain systems more intuitive and comprehensible. In order to describe the nested structure in the system we can use a tree. Each node in this tree corresponds to a nested sub-block, NSB,  in the system. In such a description the node at the root of the tree corresponds to the NSB that does not appear in the description of any other NSB. The children for each node correspond to the NSBs that appear in its description. Within this tree then the blocks at the leaves of the tree correspond to ones that does not include any other NSBs. Consequently, the description for these blocks only relies on the fundamental building blocks. We refer to this description as the fundamental description.

\begin{figure}
        \centering
        \includegraphics[width=11cm]{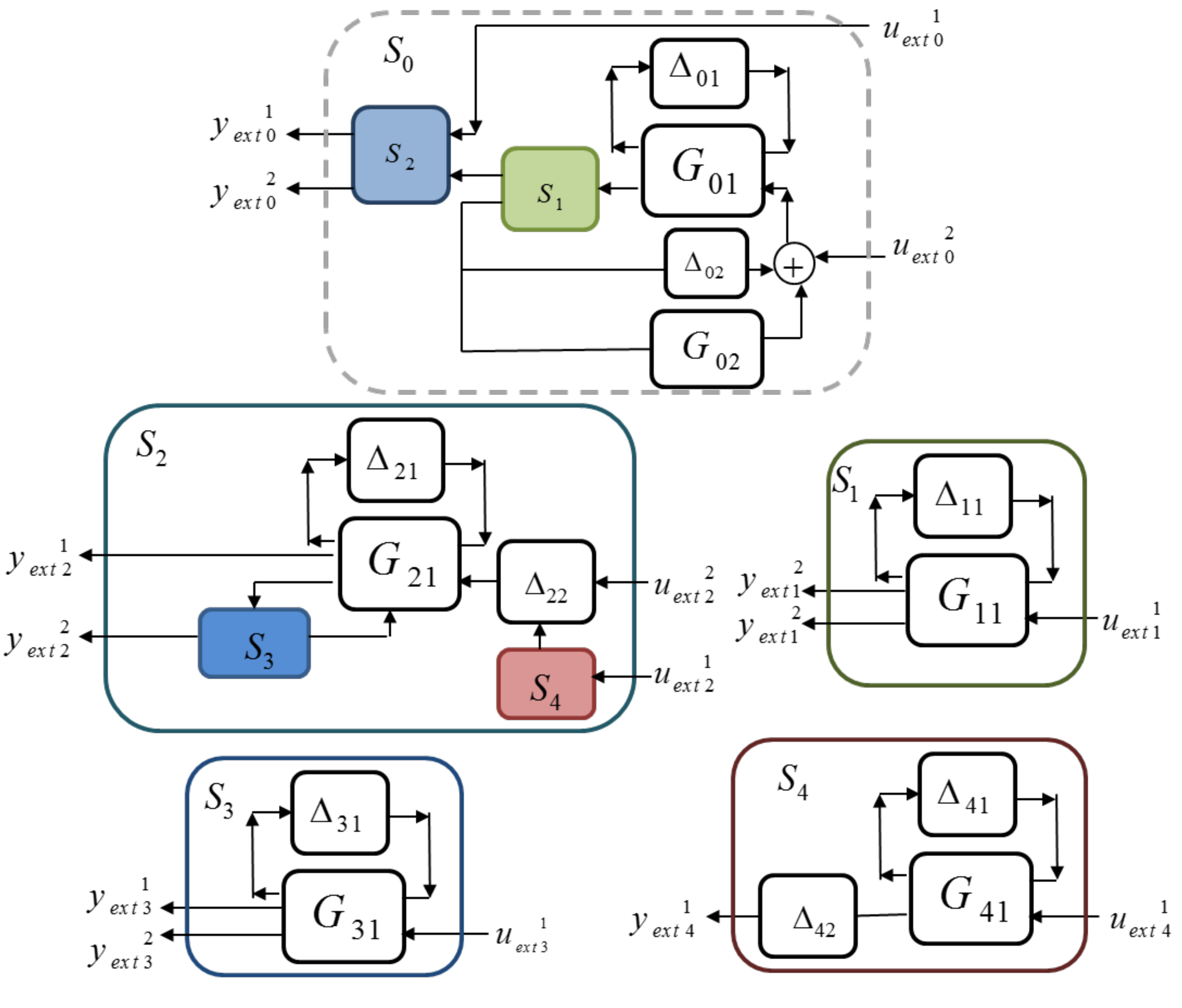}
        \caption{Blocks details of example}
        \label{fig:10}
\end{figure}

\begin{figure}
        \centering
        \includegraphics[width=2.5cm]{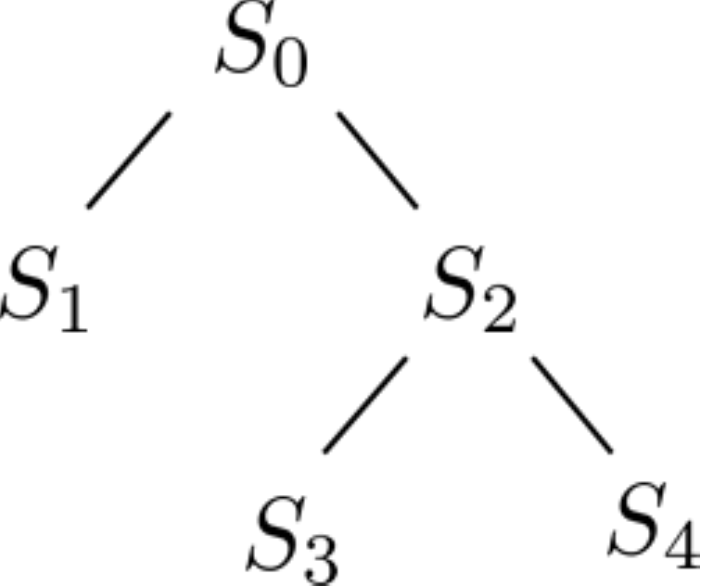}
        \caption{Example of nested structure of an interconnected system}
        \label{fig:9}
\end{figure}

At this point, let us consider an example. Consider the interconnected uncertain system depicted in Figure \ref{fig:10}. As can be seen from the figure, in this system description there are 5 NSBs, namely $S_0, S_1, \dots, S_4$. Hence the tree that describes the nested structure in this system is comprised of 5 nodes. This tree is illustrated in Figure \ref{fig:9}. As can be seen from the figure, since $S_0$ does not appear in the description of any other NSBs it is assigned to the root of the tree. Also in this tree, for instance, $S_3$ is child of $S_2$ since $S_3$ appears in the description of $S_2$.

This tree-based description of the nested structure in the system is at the heart of our approach in retrieving the interconnection matrix of the system. Before we present our algorithm for this purpose, we first need to provide a formal description for each of the NSBs. Let us consider NSB $i$ and assume that the dynamic system blocks and uncertainty blocks in its description are given by $G^i$ and $\Delta^i$ which are both block diagonal and each diagonal entry corresponds to each of such blocks in this NSB. Also let us denote the inputs and outputs vectors from the NSBs within $i$th NSB description as $u^i_s$ and $y^i_s$, i.e., these vectors correspond to the inputs and outputs from the underlying NSBs stacked. We can then summarize the mathematical description of the $i$th NSB as
\begin{subequations} \label{EQ_14abc_}
\begin{align}
y^i&=G^iu^i,\\
q^i&=\Delta^ip^i,\\
\bar u^i &= \bar \Gamma^i \bar y^i.
\end{align}
\end{subequations}
where $\bar u^i = (u^i, p^i, u^i_s, y^i_{\textrm{ext}})$ and $\bar y^i = (y^i, q^i, y^i_s, u^i_{\textrm{ext}})$ with $u^i_{\textrm{ext}}$ and $ y^i_{\textrm{ext}}$ denoting the external input and output for the $i$th NSB. Similar to \eqref {EQ_8_} the matrix $\bar \Gamma^i$ can be given as
\begin{align}\label{EQ_15_} \centering
\bar\Gamma^i = \begin{bmatrix} \Gamma_{gg}^i& \Gamma_{gd}^i& \Gamma_{gs}^i& \Gamma_{gi}^i\\ \Gamma_{dg} ^i& \Gamma_{dd}^i & \Gamma_{ds}^i & \Gamma_{di}^i\\ \Gamma_{sg}^i & \Gamma_{sd}^i & \Gamma_{ss}^i & \Gamma_{si}^i \\ \Gamma_{og}^i & \Gamma_{od}^i & \Gamma_{os}^i & \Gamma_{oi}^i\end{bmatrix}.
\end{align}
Here we assume that $\Gamma_{oi}^i = 0$. This assumption implies that there is no direct link between the external inputs and outputs. Notice that this assumption is not restrictive and can be imposed by simply adding dummy blocks and input-outputs. As was mentioned before, the description for each NSB can be expressed in terms of fundamental building blocks. In order to see this, notice that the description of NSBs at the leaves of the tree are given in terms of fundamental blocks. This means that the fundamental description for these blocks are readily available. For instance for the example in Figure~\ref{fig:10} the description of NSB $S_3$, that is a leaf of the tree, is given as
\begin{subequations}
\begin{align}
y^{3}&=G_{31}u^{3},\\
q^{3}&=\Delta_{31}p^{3},\\
\bar u^{3} &= \bar{\mathbf{\Gamma}}^{3}  \bar y^{3}.
\end{align}
\end{subequations}
which is in the same form as \eqref{EQ_14abc_} and only depends on fundamental building blocks. The only difference from \eqref{EQ_14abc_} is that in $\bar u^{3}$, $\bar y^{3}$ the inputs and outputs concerning the underlying NSBs are missing as they do not exist. Particularly, $\bar u^{3} = (u^{3}, p^{3}, y^{3}_{\textrm{ext}})$ and $\bar y^{3} = (y^{3}, q^{3}, u^{3}_{\textrm{ext}})$, and consequently, based on the description of $S_3$ in Figure \ref{fig:10},
\begin{align}
\bar{\mathbf \Gamma}^3 = \begin{bmatrix}\begin{array}{ccc:c:c} 0 & 0 & 0 & 1 & 0 \\ 0 & 0 & 0 & 0 & 1 \\ \hdashline 1 & 0 & 0 & 0 & 0 \\ \hdashline 0 & 1 & 0 & 0 & 0 \\ 0 & 0 & 1 & 0 & 0   \end{array}  \end{bmatrix}.
\end{align}
Having described the NSBs at the leaves of the tree, it is then possible for each parent to the leaves to merge their description with that of its own. For instance let an NSB $k$ be a parent to only leaf nodes, indexed as $i_1, \dots, i_k$ and with descriptions given as in \eqref{EQ_14abc_}. The $k$th NSB can then merge the description of its children (the leaves) NSBs with that of its own as
\begin{subequations}\label{eq:merging}
\begin{align}
\begin{bmatrix} y^k \\ y^{i_1}\\ y^{i_2} \\ \vdots \\ y^{i_k} \end{bmatrix} &= \blkdiag(G^k, G^{i_1}, G^{i_2}, \dots, G^{i_k}) \begin{bmatrix}u^k \\ u^{i_1}\\ u^{i_2} \\ \vdots \\ u^{i_k} \end{bmatrix}\\
\begin{bmatrix} q^k \\ q^{i_1}\\ q^{i_2} \\ \vdots \\ q^{i_k} \end{bmatrix} &= \blkdiag(\Delta^k, \Delta^{i_1}, \Delta^{i_2}, \dots, \Delta^{i_k}) \begin{bmatrix}p^k \\ p^{i_1}\\ p^{i_2} \\ \vdots \\ p^{i_k} \end{bmatrix}\\
\bar u^{k} = \bar \Gamma^{k} &\bar y^{k}, \quad \bar u^{i_1} = \bar{\mathbf \Gamma}^{i_1} \bar y^{i_k}, \quad \dots, \quad \bar u^{i_k} = \bar{\mathbf \Gamma}^{i_k} \bar y^{i_k} \label{eq:merging-c}
\end{align}
\end{subequations}
At this point, it is also possible to write the description of NSB $k$ only in terms of the fundamental building blocks as for the leaves. To this end, we need to eliminate the external inputs and outputs for NSBs $i_1, \dots, i_k$. This can be done using the interconnection descriptions in \eqref{eq:merging-c}. Let us illustrate this using the example in Figure \ref{fig:10}. As for NSB $S_3$, we can also describe NSB $S_4$ as
\begin{subequations}
\begin{align}
y^{4}&=G_{41}u^{4},\\
q^{4}&= \begin{bmatrix} \Delta_{41} & 0\\0 & \Delta_{42} \end{bmatrix}p^{4},\\
\bar u^{4} &= \bar{\mathbf \Gamma}^{4} \bar y^{4},
\end{align}
\end{subequations}
with
\begin{align}
\bar{\mathbf \Gamma}^4 = \begin{bmatrix} \mathbf\Gamma_{gg}^4 & \mathbf\Gamma^4_{gd} & \mathbf\Gamma^4_{gi} \\ \mathbf\Gamma^4_{dg} & \mathbf\Gamma^4_{dd} & \mathbf\Gamma^4_{di} \\ \mathbf\Gamma^4_{og} & \mathbf\Gamma^4_{od} & \mathbf\Gamma^4_{oi}  \end{bmatrix} := \begin{bmatrix}\begin{array}{cc:cc:c} 0 & 0 & 1 & 0 & 0 \\ 0 & 0 & 0 & 0 & 1 \\ \hdashline 1 & 0 & 0 & 0 & 0 \\  0 & 1 & 0 & 0 & 0 \\ \hdashline 0 & 0 & 0 & 1 & 0   \end{array}  \end{bmatrix}.
\end{align}
At this point we can describe NSB $S_2$ as in \eqref{eq:merging} with
\begin{subequations}\label{eq:mergeExample}
\begin{align}
\begin{bmatrix} y^2 \\ y^{3}\\ y^{4}  \end{bmatrix} &= \blkdiag \left(G^2, G^{3}, G^{4}\right) \begin{bmatrix}u^2 \\ u^{3}\\ u^{4}\end{bmatrix}\\
\begin{bmatrix} q^2 \\ q^{3}\\ q^{4} \end{bmatrix} &= \blkdiag(\Delta^2, \Delta^{3}, \Delta^{4}) \begin{bmatrix}p^2 \\ p^{3}\\ p^{4} \end{bmatrix}\\
\bar u^{2} = \bar \Gamma^{2} &\bar y^{2}, \quad \bar u^{3} = \bar{\mathbf \Gamma}^{3} \bar y^{3}, \quad \bar u^{4} = \bar{\mathbf \Gamma}^{4} \bar y^{4} \label{eq:mergeExample-c}
\end{align}
\end{subequations}
where $G^2 = G_{21}$ and $\Delta^2 = \blkdiag(\Delta_{21}, \Delta_{22})$. Let us now stack the interconnection constraints in \eqref{eq:mergeExample-c} as
\begin{align}
\begin{bmatrix} \begin{array}{c} u^2\\ p^2 \\ u^2_s\\ y_{\textrm{ext}}^2 \\ \hdashline u^3\\ p^3 \\ y_{\textrm{ext}}^3 \\ \hdashline u^4\\ p^4 \\ y_{\textrm{ext}}^4  \end{array} \end{bmatrix} = \begin{bmatrix} \begin{array}{c:c:c} \bar\Gamma^2 & 0 & 0\\ \hdashline 0 &\bar{\mathbf \Gamma}^3 & 0\\ \hdashline 0& 0& \bar{\mathbf \Gamma}^4 \end{array} \end{bmatrix}\begin{bmatrix} \begin{array}{c} y^2\\ q^2 \\ y^2_s\\ u_{\textrm{ext}}^2 \\ \hdashline y^3\\ q^3 \\ u_{\textrm{ext}}^3 \\ \hdashline y^4\\ q^4 \\ u_{\textrm{ext}}^4  \end{array} \end{bmatrix}.
\end{align}
Let us reorder blocks in this equation and rewrite it as
\begin{align}\label{eq:TotalInter}
\begin{bmatrix} \begin{array}{c} u^2\\ u^3 \\ u^4 \\ \hdashline p^2 \\ p^3 \\ p^4 \\ \hdashline  y_{\textrm{ext}}^2 \\ \hdashline u^2_s \\ \hdashline  y_{\textrm{ext}}^3 \\  y_{\textrm{ext}}^4  \end{array} \end{bmatrix}
= \begin{bmatrix} \begin{array}{ccc:ccc:c:c:cc} \Gamma_{gg}^2 & 0 & 0 & \Gamma_{gd}^2 & 0 & 0 & \Gamma_{gi}^2 & \Gamma_{gs}^2 & 0 & 0 \\ 0 & \mathbf\Gamma_{gg}^3 & 0 & 0 &\mathbf\Gamma_{gd}^3 & 0 & 0 & 0 & \mathbf\Gamma_{gi}^3 & 0 \\ 0 & 0 & \mathbf\Gamma_{gg}^4 & 0 & 0 &\mathbf\Gamma_{gd}^4 & 0 & 0 & 0 & \mathbf\Gamma_{gi}^4 \\ \hdashline \Gamma_{dg}^2 & 0 & 0 & \Gamma_{dd}^2 & 0 & 0 & \Gamma_{di}^2 & \Gamma_{ds}^2 & 0 & 0 \\ 0 & \mathbf\Gamma_{dg}^3 & 0 & 0 &\mathbf\Gamma_{dd}^3 & 0 & 0 & 0 & \mathbf\Gamma_{di}^3 & 0 \\ 0 & 0 & \mathbf\Gamma_{dg}^4 & 0 & 0 &\mathbf\Gamma_{dd}^4 & 0 & 0 & 0 & \mathbf\Gamma_{di}^4 \\ \hdashline \Gamma_{og}^2 & 0 & 0 & \Gamma_{od}^2 & 0 & 0 & 0 & \Gamma_{os}^2 & 0 & 0 \\ \hdashline \Gamma_{sg}^2 & 0 & 0 & \Gamma_{sd}^2 & 0 & 0 & \Gamma^2_{si} & \Gamma_{ss}^2 & 0 & 0 \\ \hdashline 0 & \mathbf\Gamma_{og}^3 & 0 & 0 &\mathbf\Gamma_{od}^3 & 0 & 0 & 0 & 0 & 0 \\ 0 & 0 & \mathbf\Gamma_{og}^4 & 0 & 0 &\mathbf\Gamma_{od}^4 & 0 & 0 & 0 & 0 \end{array} \end{bmatrix}
\begin{bmatrix} \begin{array}{c} y^2\\ y^3 \\ y^4 \\ \hdashline q^2 \\ q^3 \\ q^4 \\ \hdashline  u_{\textrm{ext}}^2 \\ \hdashline y^2_s \\ \hdashline  u_{\textrm{ext}}^3 \\  u_{\textrm{ext}}^4  \end{array} \end{bmatrix}.
\end{align}
Notice that in order to convert this description to a fundamental one, we need to eliminate $y_{\textrm{ext}}^3, y_{\textrm{ext}}^4, u_{\textrm{ext}}^3, u_{\textrm{ext}}^4, u_s^2, y_s^2$ or the last three block equations from this description. This can be done by firstly noting that
\begin{align}\label{eq:merge1}
u_s^2 = \begin{bmatrix} u_{\textrm{ext}}^3 \\ u_{\textrm{ext}}^4 \end{bmatrix}, \quad y_s^2 = \begin{bmatrix} y_{\textrm{ext}}^3 \\ y_{\textrm{ext}}^4 \end{bmatrix}.
\end{align}
Recall that we have
\begin{align}
u_s^2 = \Gamma_{sg}^2y^2 + \Gamma_{sd}^2q^2 + \Gamma_{ss}^2y_s^2 + \Gamma_{si}^2u_{\textrm{ext}}^2,
\end{align}
and
\begin{align}\label{eq:merge2}
\begin{bmatrix} y_{\textrm{ext}}^3 \\ y_{\textrm{ext}}^4 \end{bmatrix} &= \begin{bmatrix} \mathbf\Gamma^3_{og} & \mathbf\Gamma^3_{od} & \mathbf\Gamma^3_{oi} & 0 & 0 & 0 \\ 0 & 0 & 0 & \mathbf\Gamma^4_{og} & \mathbf\Gamma^{4}_{od} & \mathbf\Gamma^4_{oi}\end{bmatrix} \begin{bmatrix} y^3 \\ q^3 \\ u_{\textrm{ext}}^3 \\ y^4 \\ q^4 \\ u_{\textrm{ext}}^4 \end{bmatrix}\notag\\
&= \begin{bmatrix} \mathbf\Gamma^3_{og} & 0 \\ 0 & \mathbf\Gamma^4_{og}\end{bmatrix} \begin{bmatrix} y^3  \\ y^4 \end{bmatrix} + \begin{bmatrix} \mathbf\Gamma^3_{od} & 0\\ 0&  \mathbf\Gamma^{4}_{od}\end{bmatrix} \begin{bmatrix}  q^3 \\ q^4\end{bmatrix}.
\end{align}
Combining equations in \eqref{eq:merge1}--\eqref{eq:merge2} will result in
\begin{equation}
\begin{split}
\begin{bmatrix} u_{\textrm{ext}}^3 \\ u_{\textrm{ext}}^4 \end{bmatrix} &= \begin{bmatrix}  \Gamma^2_{sg} &  \Gamma_{sd}^2 &  \Gamma^2_{si} \end{bmatrix} \begin{bmatrix} y^2\\ q^2\\ u^2_{\textrm{ext}} \end{bmatrix} + \Gamma_{ss}^2 \begin{bmatrix} y_{\textrm{ext}}^3 \\ y_{\textrm{ext}}^4 \end{bmatrix} \\ & = \begin{bmatrix}  \Gamma^2_{sg} &  \Gamma_{sd}^2 &  \Gamma^2_{si} \end{bmatrix} \begin{bmatrix} y^2\\ q^2\\ u^2_{\textrm{ext}} \end{bmatrix} + \Gamma_{ss}^2\left( \begin{bmatrix}\mathbf \Gamma^3_{og} & 0 \\ 0 &\mathbf \Gamma^4_{og}\end{bmatrix} \begin{bmatrix} y^3  \\ y^4 \end{bmatrix} + \begin{bmatrix}\mathbf \Gamma^3_{od} & 0\\ 0& \mathbf \Gamma^{4}_{od}\end{bmatrix} \begin{bmatrix}  q^3 \\ q^4\end{bmatrix} \right)\\
& =: \Gamma_{sgdi}^2 \begin{bmatrix} y^2\\ q^2\\ u^2_{\textrm{ext}} \end{bmatrix} \Gamma_{ss}^2\left( \Gamma_{ogc}^2 y_c^2 + \Gamma_{odc}^2q_c^2 \right)
\end{split}
\end{equation}
This then allows us to rewrite the set of equations in \eqref{eq:TotalInter} that we want to keep as
\begin{equation}\label{eq:mergingexample1}
\begin{split}
\begin{bmatrix} u^3 \\ u^4 \end{bmatrix} =& \left( \begin{bmatrix}\mathbf \Gamma_{gg}^3 & 0 \\ 0 &\mathbf \Gamma_{gg}^4 \end{bmatrix} + \begin{bmatrix} \mathbf\Gamma_{gi}^3 & 0 \\ 0 & \mathbf\Gamma_{gi}^4 \end{bmatrix}\Gamma_{ss}^2\begin{bmatrix} \mathbf\Gamma_{og}^3 & 0 \\ 0 & \mathbf\Gamma_{og}^4 \end{bmatrix} \right) \begin{bmatrix} y^3 \\ y^4 \end{bmatrix} + \\& \left( \begin{bmatrix} \mathbf\Gamma_{gd}^3 & 0 \\ 0 & \mathbf\Gamma_{gd}^4 \end{bmatrix} + \begin{bmatrix} \mathbf\Gamma_{gi}^3 & 0 \\ 0 & \mathbf\Gamma_{gi}^4 \end{bmatrix}\Gamma_{ss}^2\begin{bmatrix} \mathbf\Gamma_{od}^3 & 0 \\ 0 & \mathbf\Gamma_{od}^4 \end{bmatrix} \right) \begin{bmatrix} q^3 \\ q^4 \end{bmatrix}+\\& \begin{bmatrix} \mathbf\Gamma_{gi}^3 & 0 \\ 0 & \mathbf\Gamma_{gi}^4 \end{bmatrix} \begin{bmatrix}  \Gamma^2_{sg} &  \Gamma_{sd}^2 &  \Gamma^2_{si} \end{bmatrix} \begin{bmatrix} y^2\\ q^2\\ u^2_{\textrm{ext}}\end{bmatrix} \\ =&: \left( \Gamma_{ggc}^2 + \Gamma_{gic}^2 \Gamma_{ss}^2 \Gamma_{ogc}^2 \right)y_c + \left( \Gamma_{gdc}^2 + \Gamma_{gic}^2 \Gamma_{ss}^2 \Gamma_{odc}^2 \right)q_c + \Gamma_{gic}^2\Gamma_{sgdi}^2 \begin{bmatrix} y^2\\ q^2\\ u^2_{\textrm{ext}} \end{bmatrix}
\end{split}
\end{equation}
\begin{equation}
\begin{split}
u^2 &= \begin{bmatrix}  \Gamma^2_{gg} &  \Gamma_{gd}^2 &  \Gamma^2_{gi} \end{bmatrix} \begin{bmatrix} y^2\\ q^2\\ u^2_{\textrm{ext}} \end{bmatrix} + \Gamma_{gs}^2\left( \Gamma_{ogc}^2 \begin{bmatrix} y^3  \\ y^4 \end{bmatrix} + \Gamma_{odc}^2 \begin{bmatrix}  q^3 \\ q^4\end{bmatrix} \right)\\
& =: \Gamma_{ggdi}^2 \begin{bmatrix} y^2\\ q^2\\ u^2_{\textrm{ext}} \end{bmatrix} + \Gamma_{gs}^2\left( \Gamma_{ogc}^2 y_c^2 + \Gamma_{odc}^2 q_c^2 \right)
\end{split}
\end{equation}
\begin{equation}
\begin{split}
\begin{bmatrix} p^3 \\ p^4 \end{bmatrix} =& \left( \begin{bmatrix}\mathbf \Gamma_{dg}^3 & 0 \\ 0 &\mathbf \Gamma_{dg}^4 \end{bmatrix} + \begin{bmatrix} \mathbf\Gamma_{di}^3 & 0 \\ 0 & \mathbf\Gamma_{di}^4 \end{bmatrix}\Gamma_{ss}^2\begin{bmatrix}\mathbf \Gamma_{di}^3 & 0 \\ 0 & \mathbf\Gamma_{di}^4 \end{bmatrix} \right) \begin{bmatrix} y^3 \\ y^4 \end{bmatrix} + \\& \left( \begin{bmatrix} \mathbf\Gamma_{dd}^3 & 0 \\ 0 & \mathbf\Gamma_{dd}^4 \end{bmatrix} + \begin{bmatrix} \mathbf\Gamma_{di}^3 & 0 \\ 0 & \mathbf\Gamma_{di}^4 \end{bmatrix}\Gamma_{ss}^2\begin{bmatrix} \mathbf\Gamma_{od}^3 & 0 \\ 0 & \mathbf\Gamma_{od}^4 \end{bmatrix} \right) \begin{bmatrix} q^3 \\ q^4 \end{bmatrix}+\\& \begin{bmatrix} \mathbf\Gamma_{di}^3 & 0 \\ 0 & \mathbf\Gamma_{di}^4 \end{bmatrix} \begin{bmatrix}  \Gamma^2_{sg} &  \Gamma_{sd}^2 &  \Gamma^2_{si} \end{bmatrix} \begin{bmatrix} y^2\\ q^2\\ u^2_{\textrm{ext}}\end{bmatrix} \\ =&: \left( \Gamma_{dgc}^2 + \Gamma_{dic}^2 \Gamma_{ss}^2 \Gamma_{ogc}^2 \right)y_c + \left( \Gamma_{ddc}^2 + \Gamma_{dic}^2 \Gamma_{ss}^2 \Gamma_{odc}^2 \right)q_c + \Gamma_{dic}^2\Gamma_{sgdi}^2 \begin{bmatrix} y^2\\ q^2\\ u^2_{\textrm{ext}} \end{bmatrix}
\end{split}
\end{equation}
\begin{equation}
\begin{split}
p^2 &= \begin{bmatrix}  \Gamma^2_{dg} &  \Gamma_{dd}^2 &  \Gamma^2_{di} \end{bmatrix} \begin{bmatrix} y^2\\ q^2\\ u^2_{\textrm{ext}} \end{bmatrix} + \Gamma_{ds}^2\left( \Gamma_{ogc}^2 \begin{bmatrix} y^3  \\ y^4 \end{bmatrix} + \Gamma_{odc}^2 \begin{bmatrix}  q^3 \\ q^4\end{bmatrix} \right)\\
& =: \Gamma_{dgdi}^2 \begin{bmatrix} y^2\\ q^2\\ u^2_{\textrm{ext}} \end{bmatrix} + \Gamma_{ds}^2\left( \Gamma_{ogc}^2 y_c^2 + \Gamma_{odc}^2 q_c^2 \right)
\end{split}
\end{equation}
and
\begin{equation}\label{eq:mergingexample2}
\begin{split}
y_{\textrm{ext}}^2 &= \begin{bmatrix}  \Gamma^2_{og} &  \Gamma_{od}^2 \end{bmatrix} \begin{bmatrix} y^2\\ q^2\end{bmatrix} + \Gamma_{os}^2\left( \Gamma_{ogc}^2 \begin{bmatrix} y^3  \\ y^4 \end{bmatrix} + \Gamma_{odc}^2 \begin{bmatrix}  q^3 \\ q^4\end{bmatrix} \right)\\
& =: \Gamma_{ogd}^2 \begin{bmatrix} y^2\\ q^2\end{bmatrix} + \Gamma_{os}^2\left( \Gamma_{ogc}^2 y_c^2 + \Gamma_{odc}^2 q_c^2 \right)
\end{split}
\end{equation}
Notice that this description is in the same format as the description given for the leaves and in fact can be written as
\begin{subequations}
\begin{align}
\begin{bmatrix} u^2\\ u^3 \\ u^4 \end{bmatrix} &= \diag(G^2, G^3, G^4) \begin{bmatrix} y^2\\ y^3 \\ y^4 \end{bmatrix},\\
\begin{bmatrix} q^2\\ q^3 \\ q^4 \end{bmatrix} &= \diag(\Delta^2, \Delta^3, \Delta^4) \begin{bmatrix} p^2\\ p^3 \\ p^4 \end{bmatrix},\\
\begin{bmatrix} u^2\\ u^3 \\ u^4\\ p^2\\ p^3 \\ p^4\\ y_{\textrm{ext}}^2 \end{bmatrix} &= \bar{\mathbf\Gamma}^i \begin{bmatrix} y^2\\ y^3 \\ y^4 \\ q^2\\ q^3 \\ q^4 \\ u_{\textrm{ext}}^2 \end{bmatrix}
\end{align}
\end{subequations}
by partitioning the defined matrices in \eqref{eq:mergingexample1}--\eqref{eq:mergingexample2} accordingly, which is in fundamental form. The same procedure can be generalized for all problems, that is for an NSB $k$ with all its children descriptions given in fundamental form, it is possible to recover the fundamental description for this NSB by simply forming
\small
\begin{align*}
\Gamma_{sgdi}^k &= \begin{bmatrix}  \Gamma^k_{sg} &  \Gamma_{sd}^k &  \Gamma^k_{si} \end{bmatrix}, \quad \Gamma_{ggdi}^k = \begin{bmatrix}  \Gamma^k_{gg} &  \Gamma_{gd}^k &  \Gamma^k_{gi} \end{bmatrix}, \quad \Gamma_{dgdi}^k = \begin{bmatrix}  \Gamma^k_{dg} &  \Gamma_{dd}^k &  \Gamma^k_{di} \end{bmatrix}\\
\Gamma_{ggc}^k &= \diag(\mathbf\Gamma_{gg}^{i_1}, \dots,\mathbf\Gamma_{gg}^{i_k}), \quad \Gamma_{gic}^k = \diag(\mathbf\Gamma_{gi}^{i_1}, \dots,\mathbf\Gamma_{gi}^{i_k}),\quad \Gamma_{ogc}^k = \diag(\mathbf\Gamma_{og}^{i_1}, \dots,\mathbf\Gamma_{og}^{i_k})\\
\Gamma_{gdc}^k &= \diag(\mathbf\Gamma_{gd}^{i_1}, \dots,\mathbf\Gamma_{gd}^{i_k}), \quad \Gamma_{odc}^k = \diag(\mathbf\Gamma_{od}^{i_1}, \dots,\mathbf\Gamma_{od}^{i_k}),\quad \Gamma_{dgc}^k = \diag(\mathbf\Gamma_{dg}^{i_1}, \dots,\mathbf\Gamma_{dg}^{i_k})\\
\Gamma_{dic}^k &= \diag(\mathbf\Gamma_{di}^{i_1}, \dots,\mathbf\Gamma_{di}^{i_k}), \quad \Gamma_{ddc}^k = \diag(\mathbf\Gamma_{dd}^{i_1}, \dots,\mathbf\Gamma_{dd}^{i_k}),
\end{align*}
\normalsize
where $\{ i_1, \dots, i_k \}$ are the children of NSB $k$, and through equations \eqref{eq:mergingexample1}--\eqref{eq:mergingexample2}. Following the same procedure upwards the tree, the fundamental description at the root will correspond to the description of the interconnected system as in \eqref{EQ_10abc_}.

\section{Conclusions}
\label{sec:4}
In this paper we put forth an automated framework for extracting standard description of large-scale interconnected uncertain systems used in scalable analysis algorithms. Particularly we showed that given a nested block-diagram-based description of an interconnected uncertain system, it is possible to first describe the nested structure in the system using a tree, and how based on this it is possible to recover a suitable description of the system.

\end{document}